\documentclass[pdflatex,sn-mathphys-num]{sn-jnl}

\usepackage{graphicx}%
\usepackage{multirow}%
\usepackage{amsmath,amssymb,amsfonts}%
\usepackage{amsthm}%
\usepackage{mathrsfs}%
\usepackage[title]{appendix}%
\usepackage{xcolor}%
\usepackage{textcomp}%
\usepackage{manyfoot}%
\usepackage{booktabs}%
\usepackage{algorithm}%
\usepackage{algorithmicx}%
\usepackage{algpseudocode}%
\usepackage{listings}%

\theoremstyle{thmstyleone}%
\newtheorem{theorem}{Theorem}
\theoremstyle{thmstyletwo}%

\theoremstyle{thmstylethree}%
\newtheorem{definition}{Definition}%

\newcommand{\up}[1]{^{\mathrm{#1}}}
\usepackage{eurosym}
\usepackage{tikz}
\usepackage{tikzscale}
\usetikzlibrary{patterns}
\usetikzlibrary{datavisualization.formats.functions}
\usetikzlibrary{decorations.pathreplacing,calligraphy}
\usetikzlibrary{backgrounds, positioning, shapes.geometric, calc}
\usepackage[caption=false,font=normalsize,labelfont=sf,textfont=sf]{subfig}
\usepackage{pdflscape}

\begin{document}

\title[Virtual Linking Bids for Market Clearing with Non-Merchant Storage]{Virtual Linking Bids for Market Clearing with Non-Merchant Storage}


\author*[1]{\fnm{Eléa} \sur{Prat}}\email{emapr@dtu.dk}

\author[1]{\fnm{Jonas} \sur{Bodulv Broge}}

\author[1]{\fnm{Richard M.} \sur{Lusby}}\email{rmlu@dtu.dk}

\affil*[1]{\orgdiv{Department of Technology, Management and Economics}, \orgname{Technical University of Denmark}, \orgaddress{\city{Kgs. Lyngby}, \postcode{2800}, \country{Denmark}}}


\abstract{In the context of energy market clearing, non-merchant assets are assets that do not submit bids but 
whose operational constraints are included. 
Integrating energy storage systems as non-merchant assets can maximize social welfare. However, the disconnection between consecutive market clearings poses challenges for market properties, and this is not well studied yet.
We contribute to the literature on market clearing with non-merchant storage by proposing a market-clearing procedure that preserves desirable market properties, even under uncertainty. 
This approach is based on a novel representation of storage systems in which the energy available is discretized to reflect the different prices at which the storage system was charged. These prices are then included as virtual bids, establishing a link between different market clearings. We show that market clearing with such virtual linking bids has the advantage of guaranteeing cost recovery for market participants and can outperform traditional approaches in terms of social welfare.}




\maketitle

\bmhead{Acknowledgements}

Thanks to Linde Frölke, Pierre Pinson and Jalal Kazempour for the first discussions on the concept presented in this paper. Thanks to Pierre Pinson for reviewing and providing feedback.

\section{Introduction}

In order to enable the safe operation of future energy systems with a high share of intermittent and stochastic renewable sources of energy production, the proportion of large-scale energy storage is expected to increase significantly in the coming years~\cite{IEA2022World}. A major challenge that comes with this evolution is how to best integrate storage in energy markets. This is emphasized by an order by the Federal Energy Regulatory Commission in the United States, urging system operators to implement changes to facilitate market participation of electric storage systems~\cite{FERC2018}. 
As a step towards addressing this challenge,~\cite{Singhal2020Pricing} describes three different ways in which storage systems can be included in the market clearing procedure (MCP). In the first and second options, storage systems participate similarly to conventional generators and loads, submitting price and quantity bids. We refer to this setup as \emph{merchant storage}. In the third option, the storage systems' operational constraints are included, but they do not submit bids. Instead, they are represented as a means to transfer energy between time periods, much like power lines are represented as a means to transfer energy between connected areas~\cite{Pinson2023What}. We call this \emph{non-merchant storage}.
A considerable advantage of market clearing with non-merchant storage is that it can achieve the most economically efficient outcomes and the highest social welfare~\cite{Singhal2020Pricing}, compared to market clearing with merchant storage~\cite{Sioshansi2014When, Hartwig2016Impact}. The literature on market clearing with non-merchant storage is still limited~\cite{Singhal2020Pricing, Munoz2017Financial, Huang2018Market, Frölke2024Efficiency, Jiang2023Duality, Weibelzahl2018Effects, Zhang2020Energy}, most likely because this setup requires more modifications to the current energy markets than when considering merchant storage.


The main issue to address when clearing a market with non-merchant storage is how to represent the time-linking effect of the storage system. 
In real-time markets, which include ramping products that can span over several time periods, a common approach to account for time-linking effects is market clearing with look-ahead, where only the decisions over the first time periods, termed the \emph{decision horizon}, are implemented, while the decisions over the rest of the optimization horizon are advisory~\cite{Chen2022Pricing, Zhao2020Multi, Hua2019Pricing, Cho2022Pricing}. This enables future-aware decisions regarding the state of energy at the end of the decision horizon to be made, the quality of which depends on how far into the future one considers. This advisory horizon should be recalculated from one market clearing to the next in order to guarantee that solutions over the decision horizon are not impacted by it~\cite{Cruise2019Control}. The complex problem of the length of the optimization horizon has not been studied in market-clearing contexts.
Furthermore, this approach suffers from issues with pricing, even in a deterministic setup. The work in~\cite{Hua2019Pricing} proposes two methods to reintroduce a link to previous market clearings, but does not prevent the need for non-transparent uplift payments. In~\cite{Chen2022Pricing}, an argument for non-uniform pricing is made, but the prices obtained are highly dependent on the forecasts used.
Because of these limitations, we focus on MCPs without look-ahead and with uniform pricing.
Moreover, this is the usual approach for day-ahead auctions. It is, therefore, important to study.

When clearing a market with non-merchant storage and without look-ahead, consideration of future market clearings is moved to the choice of end-of-horizon storage parameters.
In the literature on market clearing with non-merchant storage, a very common assumption is that the storage is initially empty and the state of energy at the end of the market-clearing horizon is free\footnote{This is equivalent to setting it to zero in the absence of negative prices.}~\cite{Taylor2014Financial,Jiang2023Duality,Munoz2017Financial,Huang2018Market}. This leads to myopic decision-making regarding the state of energy of the storage at the end of each market clearing. These assumptions are not always stated~\cite{Weibelzahl2018Effects}, showing that this problem is usually disregarded.
While this approach might be satisfactory for short-term storage systems, which can fill and empty in a couple of minutes to hours, it needs to be reconsidered for long-term storage systems, which can provide energy over several days and are essential tools to increase the penetration of renewable energy~\cite{Dowling2020Role}.
The question of how to choose the storage level at the end of a market clearing is similar to the question of future-aware scheduling of energy storage systems, which has been studied from the perspective of the storage system operator (see e.g.,~\cite{Cruise2019Control, Nascimento2013Optimal, Castelli2023Robust, Deane2013Derivation}). However, in the context of market clearing, it is also crucial to consider the impact of this time-linking property on pricing.

By default, the market does not remember the price at which the storage was charged and thus considers any energy available at the beginning of a market clearing to be free. 
Therefore, the resulting market price might be lower than the price the storage system paid to charge. 
While this is not a problem when the storage system is always initially empty, it becomes critical to consider this aspect when non-myopic decisions are made regarding the final state of energy.
In~\cite{Frölke2024Efficiency}, a method is proposed to reestablish the connection between different market clearings and retrieve prices that send the proper signal to the storage system. However, the method is only valid under the unrealistic assumption of perfect information. The problem of connecting current and uncertain future time periods has also been identified in the case of remuneration of the storage system with financial storage rights~\cite{Bose2019Some}, or for demand response~\cite{Werner2021Pricing}. 
A modeling approach to reflect the transfer of energy between time periods has been introduced in~\cite{Zhang2022Pricing}. However, the authors do not discuss the valuation of the charged energy in subsequent market clearings. The problem of valuation has been studied in~\cite{Werner2021Pricing} for demand response in terms of deviation from a scenario with no flexibility, which cannot always be applied in the case of a storage system.
Unlike in markets with merchant storage, where the storage operator accounts for the time-linking effects when formulating bids, in markets with non-merchant storage, it becomes the responsibility of the market operator to ensure that the storage system recovers its costs. This problem is thus specific to non-merchant storage and critical to address in this context.

We aim to help address the unanswered question: How to best design an MCP that includes non-merchant storage so as to fully exploit the potential of storage systems to increase social welfare?
We assume that the storage level at the end of each market clearing is determined independently, and we focus on the pricing problem. 
Our main contribution is a novel MCP with non-merchant storage that ensures cost recovery for all the market participants, in particular for the storage system, even when considering uncertainty. The main idea is to remember the prices at which the storage system was charged and automatically create virtual linking bids (VLBs) to reflect these values in the next market clearings.

First, in Sect.~\ref{sec:2}, we present in more detail the challenges of MCPs with non-merchant storage with the help of an illustrative example. We then introduce the MCP with VLBs in Sect.~\ref{sec:3} and show promising results in terms of market properties in Sect.~\ref{sec:4}. We conclude in Sect.~\ref{sec:5}.

\section{Market clearing with non-merchant storage} \label{sec:2}
In this section, we first present the model for the myopic version of the standard MCP with non-merchant storage. We then show how it can be extended to include information about the future instead. Finally, we show by way of an example the limitations of standard MCPs in terms of cost recovery and maximization of social welfare.

\subsection{Standard market clearing procedure}

In order to get a better understanding of the problem, we model a stylized version of the storage system using several assumptions. First, we consider that there is only one storage system in the market. We also assume that the storage system has no leakage over time, and we disregard degradation, which is also justified by the fact that there is no consensus on how to best include it \cite{Sioshansi2021Energy}.



Here, we consider one market clearing. The set $\mathcal{T}$ gathers all the time periods $t$ of the given market clearing, where each time period has a duration of $\Delta t$ (in hours). For a day-ahead market, for instance, $\mathcal{T}$ would typically correspond to one day, and $\Delta t$ to one hour. Under the given assumptions, the standard MCP with non-merchant storage can be modeled as 
\begin{subequations} \label{prob:mc_init}
\allowdisplaybreaks
\begin{align}
    \label{eq:mc_init_obj}\max_{\mathbf{x}} \quad & \Delta t \sum_{t \in \mathcal{T}}  \left ( \sum_{l \in \mathcal{L}} U_{lt} d_{lt} - \sum_{g \in \mathcal{G}} C_{gt} p_{gt}  \right )\\
    \label{eq:mc_init_energy_bal} \text{s.t.} \quad & \sum_{l \in \mathcal{L}} d_{lt} + p_t\up{C} - p_t\up{D} - \sum_{g \in \mathcal{G}} p_{gt}  = 0  \, : \, \lambda_t &&&  \forall t \in \mathcal{T}  \\
    \label{eq:stor_bal_t} & e_t = e_{t-1} + \Delta t \, \eta\up{C} p_t\up{C} - \Delta t \frac{1}{\eta\up{D}}p_t\up{D}, &&& \forall t \in \mathcal{T}\setminus \{1\} \\
    \label{eq:stor_bal_1} & e_1 = E\up{init} +  \Delta t \, \eta\up{C} p_1\up{C} - \Delta t \frac{1}{\eta\up{D}}p_1\up{D}, \\
    \label{eq:gen_bound} & 0 \leq p_{gt} \leq \overline{P}_{gt}, &&& \forall g \in \mathcal{G}, t \in \mathcal{T} \\
    \label{eq:load_bound} & 0 \leq d_{lt} \leq \overline{D}_{lt}, &&& \forall l \in \mathcal{L}, t \in \mathcal{T} \\
    \label{eq:charge_bound} & 0 \leq  p_t\up{C} \leq \overline{P}\up{C}, &&& \forall t \in \mathcal{T}\\
    \label{eq:discharge_bound} & 0 \leq  p_t\up{D} \leq \overline{P}\up{D}, &&& \forall t \in \mathcal{T}\\
    \label{eq:stor_bound} & 0 \leq e_t \leq \overline{E}, &&& \forall t \in \mathcal{T} .
\end{align}
\end{subequations} 

Here, and in the following, $\mathbf{x}$ is a vector containing all the decision variables of the model. The variables are the quantity accepted for load $l \in \mathcal{L}$, $d_{lt}$, the quantity accepted for generator $g \in \mathcal{G}$, $p_{gt}$, and for the storage system, the state of energy $e_t$, the amount charged $p_t\up{C}$, and the amount discharged $p_t\up{D}$.
The objective function~\eqref{eq:mc_init_obj} is to maximize social welfare, which corresponds to the difference between load utilities $U_{lt}$ and generation costs $C_{gt}$ for accepted offers. 
Constraint~\eqref{eq:mc_init_energy_bal} ensures the balance between generation and consumption at~$t$.
Constraints~\eqref{eq:stor_bal_t} and~\eqref{eq:stor_bal_1} update the storage level, starting from the initial level $E\up{init}$, and accounting for inefficiencies in the charging and discharging processes, represented by the parameters $\eta\up{C}$ and $\eta\up{D}$, respectively. We also introduce the round-trip efficiency $\eta = \eta\up{C}\eta\up{D}$.
The maximum generation $\overline{P}_{gt}$ and load $\overline{D}_{lt}$ are enforced by constraints~\eqref{eq:gen_bound} and~\eqref{eq:load_bound}, respectively. 
Constraint~\eqref{eq:stor_bound} sets the maximum state of energy $\overline{E}$.
The maximum charge $\overline{P}\up{C}$ and discharge $\overline{P}\up{D}$ are enforced by constraints~\eqref{eq:charge_bound} and~\eqref{eq:discharge_bound}, respectively.
The uniform market price at $t$, received by all participants, is given by $\lambda_t$, the dual variable of~\eqref{eq:mc_init_energy_bal}.

Note that this model does not enforce exclusivity of charging and discharging modes, which can be a problem when market prices are strictly negative. However, implementing this requirement introduces non-convexities, which complicates the determination of market prices. This issue, which arises when explicitly modeling storage constraints in the MCP, should be addressed, but is outside of the scope of this paper. We rather assume that the conditions are such that there always exist optimal solutions without simultaneous charge and discharge.

We make further assumptions regarding the market operation, namely that there is perfect competition and that market participants bid their true costs.

\subsection{Approaches to avoid myopic decisions}

As mentioned in~\cite{Singhal2020Pricing} and~\cite{Frölke2024Efficiency}, the model in~\eqref{prob:mc_init} is not sufficient to avoid myopic decisions regarding the state of energy at the end of the market clearing. Indeed, with this model, it is not valuable to store energy in the final time period (unless there are strictly negative prices). 
In~\cite{Singhal2020Pricing}, several options to avoid myopic decisions are listed. 

The first one is to impose a final state of energy $E\up{end}$ at the end of the market clearing, which is determined by considering information about future market clearings. This is done by adding to~\eqref{prob:mc_init} the constraint
\begin{equation}
    \label{eq:stor_end} e_{\text{T}} = E\up{end},
\end{equation}
where $t = \text{T}$ is the last time period of the market clearing.

Another option is to steer the level to a desirable value by adding a penalty term in the objective function, with a cost $S\up{end}$ in~\euro{}/MWh. This corresponds to the estimated value of each~MWh of energy kept in store at the end of the market clearing. In~\eqref{prob:mc_init}, the objective function~\eqref{eq:mc_init_obj} becomes
\begin{equation}
    \label{eq:mc_penal_obj} \max_{\mathbf{x}} \quad \Delta t  \sum_{t \in \mathcal{T}}  \left ( \sum_{l \in \mathcal{L}} U_{lt} d_{lt} - \sum_{g \in \mathcal{G}} C_{gt} p_{gt}  \right ) - S\up{end} e_{\text{T}}.
\end{equation}

\subsection{Limitations on an illustrative example} \label{sec:ex}

We show the limitations of the standard MCPs using an example.
We consider a storage system that is initially empty with $\overline{E}= 2.5$~MWh, $\overline{P}\up{C}$ = $\overline{P}\up{D} = 3.5$~MW, and $\eta\up{C} = \eta\up{D} = 0.8$. In this example, there are two market clearings of one time period of one hour each. 
This stylized setup does not reflect an actual market design but is intended solely for illustration. It enables a clear exposition of the interactions between market clearings.
The same phenomenon could be observed with market clearings of, e.g., 24 time periods, as it results from the separation between market clearings rather than from the number of time periods in each clearing.
One load and two generators participate in the market. The related parameters are listed in Table~\ref{tab:data}. Note that there is no load in the first market clearing, such that generators can only be used to charge the storage system, if necessary. The code for this example, as well as all the other examples presented in this paper, is available online at \url{https://github.com/eleaprat/MC_non_merchant_stg}.

\begin{table}
    \caption{Load and generators parameters for the illustrative example}
    \centering
    \begin{tabular}{lrrrrrrr}
        \hline
        \multicolumn{1}{c}{MC} & \multicolumn{1}{c}{$t$} & \multicolumn{1}{c}{$U_{1t}$} & \multicolumn{1}{c}{$\overline{D}_{1t}$} & \multicolumn{1}{c}{$C_{1t}$} & \multicolumn{1}{c}{$\overline{P}_{1t}$} & \multicolumn{1}{c}{$C_{2t}$} & \multicolumn{1}{c}{$\overline{P}_{2t}$} \\ 
        \hline
        1 & 1 & 12 & 0 & 5 & 2 & 10 & 2 \\
        2 & 1 & 12 & 3 & 2 & 2 & 9 & 2 \\ 
        \hline
    \end{tabular}
    \footnotetext{Prices are in~\euro{}/MWh and quantities in~MW.}
    \label{tab:data}
\end{table}

We first solve~\eqref{prob:mc_init} with~\eqref{eq:stor_end}, where $E\up{end}= 1.25$~MWh for the first market clearing and  $E\up{end}= 0$~MWh for the second market clearing\footnote{These are found to maximize the total social welfare when optimizing for the two market clearings together.}.
The resulting dispatch and market prices are shown in Table~\ref{tab:results}.
\begin{table}
    \caption{Results for the illustrative example}
    \centering
    \begin{tabular}{llrrrrrrrr}
    \hline
    & \multicolumn{1}{c}{MC} & \multicolumn{1}{c}{$t$} & \multicolumn{1}{c}{$e_{t}$ } &  \multicolumn{1}{c}{$p\up{C}_{t}$} &  \multicolumn{1}{c}{$p\up{D}_{t}$} &\multicolumn{1}{c}{$d_{1t}$} & \multicolumn{1}{c}{$p_{1t}$} & \multicolumn{1}{c}{ $p_{2t}$} & \multicolumn{1}{r}{$\lambda_t$} \\ 
    \hline
    \multirow{2}{*}{With~\eqref{eq:stor_end}} & 1 & 1 & 1.25 & 1.56 & 0 & 0 & 1.56 & 0 & 5 \\
    & 2 & 1 & 0 & 0 & 1 & 3 & 2 & 0 & {[}2,9{]} \\
    \hline
    \multirow{2}{*}{With~\eqref{eq:mc_penal_obj}} & 1 & 1 & 0 & 0 & 0 & 0 & 0 & 0 & 0 \\
    & 2 & 1 & 0 & 0 & 0 & 3 & 2 & 1 & 9 \\ 
    \hline
    \multirow{2}{*}{Myopic} & 1 & 1 & 0 & 0 & 0 & 0 & 0 & 0 & 0 \\
    & 2 & 1 & 0 & 0 & 0 & 3 & 2 & 1 & 9 \\
    \hline
    \end{tabular}
    \footnotetext{The state of energy is in~MWh, the market prices in~\euro{}/MWh and the rest of the variables in~MW.}
    \label{tab:results}
\end{table}
We can see that in market clearing~2, there is price multiplicity, where any price between 2 and 9~\euro{}/MWh is valid. It can be problematic if the final price chosen is below $\frac{5}{\eta} = 7.8125$~\euro{}/MWh 
because the storage system would then not recover its charging cost. In fact, in market clearing~2, the information that the storage system charged at 5~\euro{}/MWh is not even accessible. As a consequence, the price is chosen without accounting for it. The total social welfare is 27.19~\euro{}.

On the other hand, we can use the penalty term by solving~\eqref{prob:mc_init} with~\eqref{eq:mc_penal_obj}. Knowing about the potential price multiplicity in market clearing~2, we set the penalty $S\up{end}=2\eta=1.28$~\euro{}/MWh in the first market clearing, to make sure the storage system will then recover its costs if the lowest price of 2~\euro{}/MWh gets selected in market clearing~2. We select $S\up{end}=0$ in market clearing~2. We can see in Table~\ref{tab:results} that in this case, the storage system does not charge at all, and the resulting social welfare is reduced to 23~\euro{}. Due to this limitation, in the rest of the paper ``standard MCP'' refers to solving~\eqref{prob:mc_init} with~\eqref{eq:stor_end}.

Finally, Table~\ref{tab:results} shows that the results from solving~\eqref{prob:mc_init} with~\eqref{eq:mc_penal_obj} are the same as solving the myopic version~\eqref{prob:mc_init} alone, for this example. As formally proven in~\cite{Frölke2024Efficiency},~\eqref{prob:mc_init} ensures cost recovery but it is inefficient as it gives the storage system incentives to deviate (it can make more money with~\eqref{prob:mc_init} augmented with~\eqref{eq:stor_end}), and it does not maximize social welfare.

\section{Market clearing with virtual linking bids} \label{sec:3}

\noindent In this section, we modify the standard MCP to ensure that the cost at which the storage system charged in one market clearing is accounted for in the following market clearings. To do so, we introduce a new representation of the non-merchant storage system.

\subsection{Inter- and intra-storage}

In order to save the value at which the storage charges in view of subsequent market clearings, we introduce the concepts of \emph{net charge} and \emph{net discharge} over a market clearing. These are indicated by the difference between the initial and the final state of energy. If it is positive, it corresponds to a net charge, and if negative, it corresponds to a net discharge.

We can consider separately the exchanges of energy within the market clearing, which we associate to an \emph{intra-storage}, and the exchanges of energy with past or future market clearings, which correspond to \emph{inter-storage}. We conceptually split the storage, as illustrated in Fig.~\ref{fig:stor_val}. In the case of net charge, this quantity is saved along with the corresponding charging price in the form of different inter-storages, one for each saved price. In Fig.~\ref{fig:stor_val}, there are three inter-storages; one with price $S_1$, one with price $S_2$, and one with price $S_3$. This price is then used as a virtual bid in future market clearings until the corresponding inter-storage is emptied. In this way, we ensure that the storage will get at least what it paid for charging. Fig.~\ref{fig:stor_val} also shows the capacity for the intra-storage, which corresponds to the available capacity at the beginning of the market clearing.

\begin{figure}
    \centering
    \includegraphics{figs/tikz/stor_val.tikz}
    \caption{Storage system at the beginning of a market clearing, with previously charged quantities in red and their respective price in parentheses}
    \label{fig:stor_val}
\end{figure}

To understand how this representation of the storage system is used, two examples are shown in Fig.~\ref{fig:net}; one for net charge and one for net discharge. In the example of net charge from Fig.~\ref{fig:ch}, we see that over the market clearing, the quantity in the inter-storage is not used. At the end of the market clearing, the quantity net charged is added to the inter-storage. The example of net discharge from Fig.~\ref{fig:disch} illustrates how intra- and inter-storages can be used independently. In the two examples, the actual storage level at the end of the market clearing determines the size of the intra-storage for the next market clearing. The details of this procedure are given in the rest of this section.
\begin{figure}
    \centering
    \subfloat[\label{fig:ch}]
    {\includegraphics{figs/tikz/net_charge.tikz}}
    \\
    \subfloat[\label{fig:disch}]
    {\includegraphics{figs/tikz/net_discharge.tikz}}
    \caption{Examples of net charge (a) and net discharge (b), for a fictional market clearing of 3 time periods of one hour each. Quantities are in~MW.}
    \label{fig:net}
\end{figure}

\subsection{Overview}
Using the representation of the storage system as the combination of an intra-storage and an inter-storage, we can introduce VLBs to ensure that the price at which the storage system charged is available in subsequent market clearings. This process resembles a usual bid but is automated\footnote{Unlike block bids, which are submitted by market participants and establish a link between consecutive time periods within a single market clearing, VLBs are automatically generated and reflect the link between consecutive market clearings.}, i.e., no participant has to submit this virtual bid. This is coherent with the fact that we are considering non-merchant storage systems. 

The formulation of VLBs is part of the MCP, as shown in Fig.~\ref{fig:vlb_proc}. This figure also shows the flow of information. First, the market is cleared for $\mathcal{T}_1$, as described in Sect.~\ref{sec:mcvlb} by~\eqref{prob:mc_vos}. 
The inter-storage is then updated as explained in Sect.~\ref{sec:update}. In the case of net discharge of an inter-storage, the update is straightforward. In the case of net charge of the intra-storage, optimization problem~\eqref{prob:val_up} is solved. This procedure gives for each value $v \in \mathcal{V}$ in the inter-storage, the energy available at the beginning of the next market clearing, $E\up{init}_v$, and the price at which the storage system was charged, $S_v$. This price can be negative in principle. Set $\mathcal{V}$ corresponds to the different values saved in the storage system, similarly to what was shown in Fig.~\ref{fig:stor_val}, and it is also updated then.
The resulting new representation of the inter-storage is automatically used in the next market clearing for $\mathcal{T}_2$, along with the information submitted by the other participants.
The parameters related to the storage characteristics are not represented in the figure because we assume they are constant from one market clearing to the next. 
\begin{figure}
    \centering
    \includegraphics{figs/tikz/vlb.tikz}
    \caption{MCP with VLBs between two market clearings $\mathcal{T}_1$ and $\mathcal{T}_2$}
    \label{fig:vlb_proc}
\end{figure}

The only external information needed from the storage system for the market clearing is the final state of energy. This information would also be needed in a standard MCP. It is essential to discuss the responsibility and method for the choice of this value, as well as how it impacts market properties. However, this is outside our scope.

For simplicity, we focus on a market clearing without network constraints and with only one non-merchant storage system. Nevertheless, the procedure presented here can easily be extended to introduce these features. In the case of multiple non-merchant storage systems, the update would be performed for each storage system individually.

\subsection{Model}
\label{sec:mcvlb}

We now introduce the model for the MCP with VLBs, based on the model in~\eqref{prob:mc_init} and including the representation of the storage system with intra- and inter-storage components.
\begin{subequations} \label{prob:mc_vos}
\allowdisplaybreaks
\begin{alignat}{3}
    \label{eq:mc_obj} \max_{\mathbf{x}} \quad & \Delta t \sum_{t \in \mathcal{T}}  \left ( \sum_{l \in \mathcal{L}} U_{lt} d_{lt} - \sum_{g \in \mathcal{G}} C_{gt} p_{gt} - \sum_{v \in \mathcal{V}} S_v p_{vt}\up{D,e}  \right )\\
    \nonumber \text{s.t.} \quad & ~\eqref{eq:gen_bound} -~\eqref{eq:load_bound},\\
    \label{eq:energy_bal} & \sum_{l \in \mathcal{L}} d_{lt} + p_t\up{C,a} - p_t\up{D,a} - \sum_{g \in \mathcal{G}} p_{gt}  - \sum_{v \in \mathcal{V}} p_{vt}\up{D,e} = 0 , &&&  \forall t \in \mathcal{T} \\
    \label{eq:stor_bal_t_intra} &  e\up{a}_t = e\up{a}_{t-1} + \Delta t \, \eta\up{C} p_t\up{C,a} - \Delta t \, \frac{1}{\eta\up{D}} p_t\up{D,a} , &&& \forall t \in \mathcal{T}\setminus \{1\}\\
    \label{eq:stor_bal_1_intra} & e\up{a}_1 = \Delta t \, \eta\up{C} p_1\up{C,a}  - \Delta t \, \frac{1}{\eta\up{D}} p_1\up{D,a} , \\
    \label{eq:stor_bal_t_inter} & e\up{e}_{vt} = e\up{e}_{v,t-1} - \Delta t \, \frac{1}{\eta\up{D}} p_{vt}\up{D,e}, &&& \forall v \in \mathcal{V}, \,  t \in \mathcal{T}\setminus \{1\}\\
    \label{eq:stor_bal_1_inter} & e\up{e}_{v,1} = E\up{init}_v - \Delta t \, \frac{1}{\eta\up{D}} p_{v,1}\up{D,e}, &&& \forall v \in \mathcal{V}\\
    \label{eq:charge_bound_intra} & 0 \leq p_t\up{C,a} \leq \overline{P}\up{C}, &&&  \forall t \in \mathcal{T}\\
    \label{eq:discharge_bound_min_intra} & p_t\up{D,a} \geq 0, &&& \forall t \in \mathcal{T}\\
    \label{eq:discharge_bound_min_inter} & p_{vt}\up{D,e} \geq 0, &&& \forall v \in \mathcal{V}, \,  \forall t \in \mathcal{T}\\
    \label{eq:discharge_bound_max_all} & p_t\up{D,a} + \sum_{v \in \mathcal{V}} p_{vt}\up{D,e} \leq \overline{P}\up{D}, &&& \forall t \in \mathcal{T}\\
    \label{eq:stor_bound_all} & 0 \leq e\up{a}_t + \sum_{v \in \mathcal{V}} e\up{e}_{vt}  \leq \overline{E}, &&& \forall t \in \mathcal{T}\\
    \label{eq:stor_pos} & e\up{e}_{vt} \geq 0, &&& \forall v \in \mathcal{V}, \,  \forall t \in \mathcal{T}\\
    \label{eq:stor_intra_lower_end} & e\up{a}_\text{T} \geq 0,\\
    \label{eq:stor_end_all} & e\up{a}_\text{T} + \sum_{v \in \mathcal{V}}  e\up{e}_{v\text{T}} \geq E\up{end} .
\end{alignat}
\end{subequations} 

There are no changes regarding loads and generators, and their decision variables. There are new variables for the storage system. For the intra-storage system, $e\up{a}_t$ gives the state of energy, $p_t\up{C,a}$ the quantity charged and $p_t\up{D,a}$ the quantity discharged. For the inter-storage system, the variables are $e\up{e}_{vt}$, the state of energy for value $v$, and $p_{vt}\up{D,e}$, the quantity discharged from inter-storage with value $v$.

The objective function~\eqref{eq:mc_obj} is modified to include the VLBs from the inter-storage, with prices $S_v$. Here and in the balance constraint~\eqref{eq:energy_bal}, the inter-storage appears similarly to conventional generators. 
Constraints~\eqref{eq:stor_bal_t_intra} and~\eqref{eq:stor_bal_1_intra} update the state of energy for the intra-storage. 
Note that the intra-storage is by definition empty at the beginning of any market clearing. However, it can occasionally take negative values, which corresponds to temporarily using some of the capacity that is reserved by the inter-storage. This is necessary to ensure that the arbitrage opportunities within a given market clearing are not limited due to the intra-storage being initially empty, while there is, in practice, some energy to discharge. To limit these operations to arbitrage within this market clearing only,~\eqref{eq:stor_intra_lower_end} specifies that the final state of energy of the intra-storage cannot be negative.
Constraints~\eqref{eq:stor_bal_t_inter} and~\eqref{eq:stor_bal_1_inter} update the state of energy for each inter-storage. 
Constraint~\eqref{eq:charge_bound_intra} gives the bounds for charge, and constraints~\eqref{eq:discharge_bound_min_intra},~\eqref{eq:discharge_bound_min_inter} and~\eqref{eq:discharge_bound_max_all} together give the bounds for discharge.
Constraint~\eqref{eq:stor_bound_all} enforces the storage capacity. The state of energy from each inter-storage is positive with~\eqref{eq:stor_pos}.
Constraint~\eqref{eq:stor_end_all} is the updated version of constraint~\eqref{eq:stor_end}. It imposes a minimum final level, but the actual final level is also impacted by the VLBs in the objective function. This formulation thus combines some of the aspects of~\eqref{eq:stor_end} and~\eqref{eq:mc_penal_obj}.

There is no charge of the inter-storage during the market clearing. In the case of net charge, the inter-storage is modified outside of the market clearing, as shown in Sect.~\ref{sec:update}.

\subsection{Update of the inter-storage} \label{sec:update}

\subsubsection{Discharge of an inter-storage}
At the end of the market clearing, the state of energy in each inter-storage is updated. If the inter-storage for value $v$ is emptied in the market clearing, the corresponding index is dropped in $\mathcal{V}$.

\subsubsection{Net charge of the intra-storage}
In the case of net charge, the quantity net charged is added to the inter-storage, with the corresponding charging price divided by $\eta$. However, since the storage might have been charging at different prices during the market clearing, the question of which price to save arises. To address this, we introduce the following optimization model
\begin{subequations} \label{prob:val_up}
\allowdisplaybreaks
\begin{align}
    \label{eq:val_obj}\min_{\mathbf{p\up{C,loc}},\mathbf{p\up{C,mov}}} \quad & \Delta t \sum_{t \in \mathcal{T}} \lambda^*_{t} (p_t\up{D,a*}-p\up{C,loc}_t) \\
    \label{eq:rev_pos} \text{s.t.} \quad & \Delta t \sum_{t \in \mathcal{T}} \lambda^*_{t} (p_t\up{D,a*}-p\up{C,loc}_t) \geq 0, \\
    \label{eq:bal_loc} & \sum_{t \in \mathcal{T}} \left(\eta\up{C} p_t\up{C,loc} - \frac{1}{\eta\up{D}} p_t\up{D,a*}\right) = 0,\\
    \label{eq:total} & p_t\up{C,loc} + p_t\up{C,mov} = p_t\up{C,a*} , \, &&& \forall t \in \mathcal{T}\\
    \label{eq:mov_pos} & p_t\up{C,loc}, \, p_t\up{C,mov} \geq 0 , \, &&& \forall t \in \mathcal{T}.
\end{align}
\end{subequations}

Some of the parameters in this model, indicated with $^*$, are obtained from the solution of the market clearing.
The idea is to split the quantity charged in the intra-storage into a \emph{local} and a \emph{moved} quantity, $p_t\up{C,loc}$ and $p_t\up{C,mov}$. Then, $p_t\up{C,mov}$ is added to the inter-storage, with the market price at $t$. The total \emph{moved} quantity must correspond to the net charge, which is equivalent to setting the total \emph{local} quantity to zero (sum of total \emph{local}  charge and discharge), in~\eqref{eq:bal_loc}. 
The sum of \emph{local}  and \emph{moved} quantities has to correspond to the quantity charged in the intra-storage at a given $t$, with~\eqref{eq:total}.
Constraint~\eqref{eq:mov_pos} ensures positivity for the charged quantities.

The objective~\eqref{eq:val_obj} is to minimize the intra-storage profit over this market clearing, while still ensuring that it is positive with~\eqref{eq:rev_pos}. 
Minimizing the profit over this market clearing for the intra-storage ensures that the prices saved 
for the inter-storage are as low as possible. This is done to prioritize the storage system in future exchanges. If we maximize instead, the values saved in the inter-storage would be higher and the storage system would be less likely to be selected in future market clearings. Note that problem~\eqref{prob:val_up} does not modify the profit of the storage system over this market clearing, but only the distribution of the energy between intra- and inter-storage.




\section{Study of market properties} \label{sec:4}

\noindent We study cost recovery and social welfare for the MCP with VLBs and compare to the standard MCP,~\eqref{prob:mc_init} with~\eqref{eq:stor_end}.
We discuss the impact of imperfect foresight on these properties. Finally, we discuss a limitation of both MCPs, opening directions for future research.
We use as a benchmark an \emph{ideal MCP}, where all market clearings are optimized together. In practice, this is not realistic, as it would imply perfect foresight and the ability to solve an infinite-horizon problem.

\subsection{Cost recovery}

There is cost recovery for a market participant, if their surplus is always non-negative\footnote{It is the case for the generators and loads, which can be formally shown by formulating their individual surplus maximization problem and its dual problem, and using the strong duality theorem, as in, e.g.,~\cite{Kazempour2018Stochastic}.}.
For a storage system, however, it is not relevant to ensure that the surplus is \emph{always} non-negative. For example, the storage system could pay to charge in one market clearing, to later discharge at a higher price and make a profit in subsequent market clearings. In this example, the surplus of the storage in the first market clearing would be negative. Hence, we redefine cost recovery for a non-merchant storage system. 

\begin{definition}[Cycle and cost recovery for non-merchant storage] \label{def}
We define a \emph{cycle} as a sequence of consecutive market clearings for which the storage system is initially empty and finally returns to this same state.
We say that there is \emph{cost recovery for a non-merchant storage system} if its surplus over a cycle is always non-negative. 
\end{definition}

For a non-merchant storage system, and for a fair comparison, it only makes sense to look at the surplus over a cycle. Otherwise, we need to know the value of the stored energy to include it, which is actually the complex problem that we are trying to solve here.

\begin{theorem}
    Assuming the existence of a cycle, the MCP with VLBs, represented by~\eqref{prob:mc_vos} and update of the inter-storage as described in Sect.~\ref{sec:update}, ensures cost recovery for the non-merchant storage system.
\end{theorem}
\begin{proof}
Constraint~\eqref{eq:rev_pos} ensures that the profit of the storage system is non-negative for the quantity exchanged over the market clearing. Regarding the quantity moved to the inter-storage, the charging price multiplied by $\eta$ is saved and later used as a bid. Discharge is never imposed because of the greater or equal to constraint on the final level~\eqref{eq:stor_end_all}. This ensures that for any quantity moved to the inter-storage, the price received will be at least equal to the charging price, also accounting for inefficiencies.
\end{proof}

Note that this result is not based on the assumption of perfect foresight and is thus valid in uncertain settings.
The assumption that the storage will eventually be empty is mild. It would be challenged, for example, considering the day-ahead market, in the case of net charge of the storage in a day with very high prices that never occur again. A complete day with very high prices is unlikely to happen without being foreseen. If foreseen, this situation would illustrate poor decision-making on the final level chosen.

The example of Sect.~\ref{sec:ex} shows that cost recovery is not ensured for the storage system in the standard MCP. The storage system starts empty and finishes empty, so the two market clearings considered form a cycle. However, the surplus of the storage in this cycle can be as low as -5.8125~\euro{}.
We now consider the same example using the MCP with VLBs.
The results are shown in Table~\ref{tab:res_prices}. At the end of the first market clearing, the inter-storage is charged with a price of 5~\euro{}/MWh, and it is discharged in the next hour. The storage is now marginal and receives at least 7.8125~\euro{}/MWh, so that the minimum surplus of the storage system over this cycle is 0~\euro{}, and cost recovery stands, while social welfare is not affected.

\begin{table}
    \caption{Results for the example of Sect.~\ref{sec:ex}, for the MCP with VLBs}
    \centering
    \begin{tabular}{lr|rrrrrrrrr}
    \hline
         \multicolumn{1}{l}{MC} & \multicolumn{1}{l|}{$t$}& \multicolumn{1}{l}{$d_{1t}$} & \multicolumn{1}{l}{$p_{1t}$} & \multicolumn{1}{l}{$p_{2t}$} & \multicolumn{1}{l}{$p_t\up{C,a}$} & \multicolumn{1}{l}{$p_t\up{D,a}$} & \multicolumn{1}{l}{$p_t\up{D,e}$} & \multicolumn{1}{l}{$e\up{a}_t$} & \multicolumn{1}{l}{$e\up{e}_t$} & \multicolumn{1}{r}{$\lambda_t$} \\ \hline
         1 & 1 & 0 & 1 & 0 & 1.5625 & 0 & 0 & 1.25 & 0 & 5 \\
         2 & 1 & 3 & 2 & 0 & 0 & 0 & 1 & 0 & 0 & [7.8125,9] \\ \hline
    \end{tabular}
    \footnotetext{The states of energy are in~MWh, the market prices in~\euro{}/MWh and the rest of the variables in~MW.}
    \label{tab:res_prices}
\end{table}

\subsection{Social welfare} \label{sec:sw}

If the level of the storage system at the end of each market clearing for the MCP with VLBs is set to a value that is optimal for the ideal MCP, the storage system will follow the same trajectory.
The social welfare will be maximal, regardless of the method considered.
This is not the case if an error is made in setting the level of the storage system at the end of the market clearing, e.g., because of imperfect information.

We focus on the different scenarios that can occur then. Let $E\up{end*}$ be the ideal value for the storage system level at the end of the market clearing.
Setting $E\up{end} < E\up{end*}$ might result in discharging the storage system in a period of low prices. Instead, it could have replaced more expensive generators or fulfilled a larger load at a later time. 
Setting $E\up{end} > E\up{end*}$ might result in charging the storage system in a period of high prices, thereby increasing the production from expensive generators.
Here, we can further identify a special case, when the price paid for charging is very high and does not occur again for a long period.
In all these cases, the social welfare is reduced.
We illustrate these scenarios with simple examples.

We look at the social welfare for the MCP with VLBs, in comparison to the one obtained for the standard MCP using the same values for the final levels $E\up{end}$. We use the outcome of the ideal MCP as a benchmark.
Since social welfare is the sum of the surpluses of all the market participants, its calculation has to be carried out over a cycle, as defined in Definition~\ref{def}. 
We compare the methods on a cycle for the MCP with VLBs\footnote{A cycle for the standard MCP is not necessarily a cycle for the MCP with VLBs but the opposite is true. If $E\up{end} = 0$, the final level for the MCP with VLBs might be higher because of~\eqref{eq:stor_end_all}, while if the final level for the MCP with VLBs is equal to zero, it means that $E\up{end} = 0$ and the final level for the standard MCP is also equal to zero.}.

In the following examples, the storage has a capacity $\overline{E}= 2.5$~MWh and is initially empty. The rest of the storage system characteristics are $\overline{P}\up{C}$ = $\overline{P}\up{D} = 3.5$~MW, and $\eta\up{C} = \eta\up{D} = 0.8$. Each market clearing consists of one time period of one hour. There is one load and one generator participating in the market.

\subsubsection{$E\up{end} < E\up{end*}$}
For this example, the rest of the data is given in Table~\ref{tab:data_sw}, including the number of market clearings. The results are in Table~\ref{tab:res_sw}, and the results in column ``$e_t$'' for the standard MCP correspond to the values set for $E\up{end}$, which are also used for the MCP with VLBs.
We can see in Table~\ref{tab:res_sw} that for market clearing~2, $E\up{end} = 0$, while $E\up{end*} = 2.5$~MWh.
As a result, with the standard MCP, the use of the generator in market clearing~2 is reduced compared to the ideal MCP. However, it is increased in market clearing~3. Even so, a part of the load cannot be fulfilled. The resulting social welfare for these three market clearings is -4.625~\euro{} for the standard MCP, compared to 14.625~\euro{} for the ideal MCP.
\begin{table}
    \caption{Data for the first example of Sect.~\ref{sec:sw}}
    \centering
    \footnotesize
    \begin{tabular}{lr|rrrr} 
    \hline
    \multicolumn{1}{l}{MC} & \multicolumn{1}{l|}{$t$} & \multicolumn{1}{l}{$D_{t}$ (MW)} & \multicolumn{1}{l}{$U_{t}$ (\euro{}/MWh)} & \multicolumn{1}{l}{$P_{t}$ (MW)} & $C_{t}$ (\euro{}/MWh)\\ 
    \hline
    1 & 1 & 1 & 5 & 4 & 5 \\
    2 & 1 & 3 & 4 & 7 & 3 \\
    3 & 1 & 3 & 10 & 2 & 9 \\ 
    \hline
    \end{tabular}
    \label{tab:data_sw}
\end{table}
\begin{sidewaystable}
    \caption{Results for the first example of Sect.~\ref{sec:sw}}
    \centering
    \footnotesize
    \begin{tabular}{lr|rrrrrr|rrrrrr|rrrrrrrr}
    \hline
    \multicolumn{1}{l}{} & \multicolumn{1}{l|}{}  & \multicolumn{6}{c|}{Ideal MCP} & \multicolumn{6}{c|}{Standard MCP} & \multicolumn{8}{c}{MCP with VLBs}\\ 
    \hline
    \multicolumn{1}{l}{MC} & \multicolumn{1}{l|}{$t$} & \multicolumn{1}{c}{$d_{t}$} & \multicolumn{1}{c}{$p_{t}$} & \multicolumn{1}{c}{$p\up{C}_{t}$} & \multicolumn{1}{c}{$p\up{D}_{t}$} & $e_t$ & \multicolumn{1}{l|}{$\lambda_t$} & $d_{t}$ & $p_{t}$ & $p\up{C}_t$ & $p\up{D}_t$  & $e_t$ & $\lambda_t$ & \multicolumn{1}{l}{$d_{t}$} & \multicolumn{1}{l}{$p_{t}$} & \multicolumn{1}{l}{$p\up{C,a}_t$} & \multicolumn{1}{l}{$p\up{D,a}_t$} & \multicolumn{1}{l}{$p\up{D,e}_t$} & \multicolumn{1}{l}{$e\up{a}_t$} & \multicolumn{1}{l}{$e\up{e}_t$} & \multicolumn{1}{l}{$\lambda_t$} \\ 
    \hline
    1 & 1 & 0 & 0 & 0 & 0 & 0  & 5 & 0.875 & 4 & 3.125 & 0 & 2.5 & 5 & 0 & 3.125 & 3.125 & 0 & 0 & 2.5 & 0 & 5 \\
    2 & 1 & 3 & 6.125 & 3.125 & 0 & 2.5 & 3 & 3 & 1 & 0 & 2 & 0 & 3 & 3 & 3 & 0 & 0 & 0 & 0 & 2.5 & 3 \\
    3 & 1 & 3 & 1 & 0 & 2 & 0 & 9 & 2 & 2 & 0 & 0 & 0 & 10 & 3 & 1 & 0 & 0 & 2 & 0 & 0 & 9 \\ 
    \hline
    \end{tabular}
    \footnotetext{Here $p\up{D,e}_t=\sum_{v \in \mathcal{V}}p\up{D,e}_{vt}$ and $e\up{e}_t=\sum_{v \in \mathcal{V}}e\up{e}_{vt}$. The state of energy is in~MWh, prices in~\euro{}/MWh, and the rest is in~MW.}
    \label{tab:res_sw}
\end{sidewaystable}

For the MCP with VLBs, the update of the inter-storage is shown in Table~\ref{tab:res_upd1}.
The value of the energy in the inter-storage in market clearing~2 is $S_1 = $ 7.8125~\euro{}/MWh, which is higher than 3~\euro{}/MWh, the cost of the generator. Therefore, it does not discharge, even though $E\up{end} = 0$. In market clearing~3, the ideal strategy can be carried out. 
The resulting social welfare is 8.375~\euro{}, which is much closer to the result of the ideal MCP.
\begin{table}
    \caption{Update of the inter-storage system after each market clearing for the first example of Sect.~\ref{sec:sw}}
    \centering
    \begin{tabular}{l|rr}
    \hline
    \multicolumn{1}{l|}{} & \multicolumn{2}{c}{$v=1$} \\ 
    \hline
    \multicolumn{1}{l|}{MC} & \multicolumn{1}{l}{$E_{v}\up{init}$ (MWh)} & \multicolumn{1}{l}{$S_v$ (\euro{}/MWh)} \\ 
    \hline
    1 & 2.5 & 7.8125 \\
    2 & 2.5 & 7.8125 \\
    3 & -   & - \\ 
    \hline
    \end{tabular}
    \footnotetext{Dashes indicate depletion.}
    \label{tab:res_upd1}
\end{table}

Therefore, in this case, the social welfare can be higher for~\eqref{prob:mc_vos} than for~\eqref{prob:mc_init} with~\eqref{eq:stor_end}. Also, observe that the formulation with VLBs ensures cost recovery for the storage system.

\subsubsection{$E\up{end} > E\up{end*}$}
In the previous example, the remaining difference with the ideal MCP is due to the fact that $E\up{end} = 2.5$~MWh in market clearing~1, while $E\up{end*} = 0$. Therefore, in both the standard MCP and the MCP with VLBs, it forces the use of the generator in market clearing~1, with a cost of 5~\euro{}/MWh, instead of increasing the use of the generator in market clearing~2, which costs 3~\euro{}/MWh. 

In both cases, the impact on social welfare is the difference in prices times the charged quantity, which is $(5-3)*3.125=6.25$~\euro{}. This corresponds exactly to the difference between the social welfare for the MCP with VLBs and for the ideal MCP.
In this example, the effect of $E\up{end} > E\up{end*}$ on the loss of social welfare is the same for the standard MCP and for the MCP with VLBs.
Next, we see another example of $E\up{end} > E\up{end*}$, where the loss of social welfare is greater for the MCP with VLBs.

\subsubsection{$E\up{end} > E\up{end*}$, with very high prices}
In this case, charge is imposed in a period of high prices, which do not occur again for a long time. This prevents the use of the full storage capacity in the following market clearings for the MCP with VLBs. Then, the social welfare can be lower for~\eqref{prob:mc_vos} than for~\eqref{prob:mc_init} with~\eqref{eq:stor_end}. We show this in a second example, for which the data is given in Table~\ref{tab:data_sw_low}, including the number of market clearings. The results are in Table~\ref{tab:res_sw_low}, and the results in column ``$e_t$'' for the standard MCP correspond to the values set for $E\up{end}$, which are also used for the MCP with VLBs. The total social welfare for these six market clearings is 1315.25~\euro{} for the standard MCP, and 1261.5~\euro{} for the MCP with VLBs, compared to 1347.75~\euro{} for the ideal MCP.
\begin{table}
    \caption{Data for the second example of Sect.~\ref{sec:sw}}
    \centering
    \begin{tabular}{lr|rrrr}
    \hline
    \multicolumn{1}{l}{MC} & \multicolumn{1}{l|}{$t$} & \multicolumn{1}{l}{$D_{t}$ (MW)} & \multicolumn{1}{l}{$U_{t}$ (\euro{}/MWh)} & \multicolumn{1}{l}{$P_{t}$ (MW)} & $C_{t}$ (\euro{}/MWh)\\ 
    \hline
    1 & 1 & 10 & 35 & 15 & 20 \\
    2 & 1 & 10 & 35 & 15 & 15 \\
    3 & 1 & 10 & 35 & 15 & 1  \\ 
    4 & 1 & 10 & 35 & 15 & 15 \\ 
    5 & 1 & 10 & 35 & 15 & 1  \\
    6 & 1 & 10 & 35 & 15 & 32 \\ 
    \hline
    \end{tabular}
    \label{tab:data_sw_low}
\end{table}

\begin{sidewaystable}
    \caption{Results for the second example of Sect.~\ref{sec:sw}}
    \centering
    \begin{tabular}{lr|rrrrrr|rrrr|rrrrrr}
    \hline
    \multicolumn{1}{l}{} & \multicolumn{1}{l|}{} & \multicolumn{6}{c|}{Ideal MCP} & \multicolumn{4}{c|}{Standard MCP} & \multicolumn{6}{c}{MCP with VLBs} \\ 
    \hline
    \multicolumn{1}{l}{MC} & \multicolumn{1}{l|}{$t$} & \multicolumn{1}{c}{$d_{t}$} & \multicolumn{1}{c}{$p_{t}$} & \multicolumn{1}{c}{$p\up{C}_{t}$} & \multicolumn{1}{c}{$p\up{D}_{t}$} & \multicolumn{1}{l}{$e_t$} & \multicolumn{1}{l|}{$\lambda_t$} & $p_{t}$ & $p\up{C}_t$ & \multicolumn{1}{c}{$p\up{D}_{t}$} & $e_t$ & \multicolumn{1}{l}{$p_{t}$} & \multicolumn{1}{l}{$p\up{C,a}_t$} & \multicolumn{1}{l}{$p\up{D,a}_t$} & \multicolumn{1}{l}{$p\up{D,e}_t$} & \multicolumn{1}{l}{$e\up{a}_t$} & \multicolumn{1}{l}{$e\up{e}_t$} \\ 
    \hline
    1 & 1 & 10 & 10 & 0 & 0 & 0 & 20 & 13.125 & 3.125 & 0 & 2.5 & 13.125 & 3.125 & 0 & 0 & 2.5 & 0 \\
    2 & 1 & 10 & 10 & 0 & 0 & 0 & 15 & 8 & 0 & 2 & 0 & 10 & 0 & 0 & 0 & 0 & 2.5 \\
    3 & 1 & 10 & 13.125 & 3.125 & 0 & 2.5 & 1 & 13.125 & 3.125 & 0 & 2.5 & 10 & 0 & 0 & 0 & 0 & 2.5 \\ 
    4 & 1 & 10 & 8 & 0 & 2 & 0 & 15 & 8 & 0 & 2 & 0 & 10  & 0 & 0 & 0 & 0 & 2.5 \\ 
    5 & 1 & 10 & 13.125 & 3.125 & 0 & 2.5 & 1 & 13.125 & 3.125 & 0 & 2.5 & 10 & 0 & 0 & 0 & 0 & 2.5 \\
    6 & 1 & 10 & 8 & 0 & 2 & 0 & 32 & 8 & 0 & 2 & 0 & 8 & 0 & 0 & 2 & 0 & 0 \\ 
    \hline
    \end{tabular}
    \footnotetext{Here $p\up{D,e}_t=\sum_{v \in \mathcal{V}}p\up{D,e}_{vt}$ and $e\up{e}_t=\sum_{v \in \mathcal{V}}e\up{e}_{vt}$. The state of energy is in~MWh, prices in~\euro{}/MWh, and the rest is in~MW. Load values $d_{t}$ and market prices $\lambda_t$ are the same for all MCPs.}
    \label{tab:res_sw_low}
\end{sidewaystable}

The difference in social welfare between the standard MCP and the MCP with VLBs is due to the fact that the inter-storage does not discharge until market clearing~6. The update of the inter-storage is shown in Table~\ref{tab:res_upd2}. The prices for market clearing~2 to market clearing~5 are too low compared to $S_1=$ 31.25~\euro{}/MWh. This prevents the use of the storage system in these market clearings, while it empties and fills multiple times with the standard MCP. This decreases social welfare.
\begin{table}
    \caption{Update of the inter-storage system after each market clearing for the second example of Sect.~\ref{sec:sw}}
    \centering
    \begin{tabular}{l|rr}
    \hline
    \multicolumn{1}{l|}{} & \multicolumn{2}{c}{$v=1$} \\ 
    \hline
    \multicolumn{1}{l|}{MC} & \multicolumn{1}{l}{$E_{v}\up{init}$(MWh)} & \multicolumn{1}{l}{$S_v$ (\euro{}/MWh)} \\ 
    \hline
    1 & 2.5 & 31.25 \\
    2 & 2.5 & 31.25 \\
    3 & 2.5 & 31.25 \\
    4 & 2.5 & 31.25 \\
    5 & 2.5 & 31.25 \\
    6 & -   & - \\ 
    \hline
    \end{tabular}
    \footnotetext{Dashes indicate depletion.}
    \label{tab:res_upd2}
\end{table}

This effect can be limited by introducing a discount on the stored value. This is shown in Table~\ref{tab:res_sw_low_dis}. In this case, we reduce by 35\% the value of the energy stored in the inter-storage after each market clearing (except after the market clearing in which it was stored). The update of the inter-storage is given in Table~\ref{tab:res_upd3}. The social welfare is 1288.38~\euro{}, as the inter-storage is now also used in market clearing~4. However, when introducing such a discount, cost recovery for the storage system is not ensured anymore.

\begin{table}
    \caption{Results for the second example of Sect.~\ref{sec:sw} for the MCP with VLBs with a discount of 35\%}
    \centering
    \begin{tabular}{lr|rrrrrr} 
    \hline
    \multicolumn{1}{l}{MC} & \multicolumn{1}{l|}{$t$} & \multicolumn{1}{l}{$p_{t}$ (MW)} & \multicolumn{1}{l}{$p\up{C,a}_t$ (MW)} & \multicolumn{1}{l}{$p\up{D,a}_t$ (MW)} & \multicolumn{1}{l}{$p\up{D,e}_t$ (MW)} & \multicolumn{1}{l}{$e\up{a}_t$ (MWh)} & \multicolumn{1}{l}{$e\up{e}_t$ (MWh)}\\ 
    \hline
    1 & 1 & 13.125 & 3.125 & 0 & 0 & 2.5 & 0 \\
    2 & 1 & 10 & 0 & 0 & 0 & 0 & 2.5 \\
    3 & 1 & 10 & 0 & 0 & 0 & 0 & 2.5 \\ 
    4 & 1 & 8 & 0 & 0 & 2 & 0 & 0 \\ 
    5 & 1 & 13.125 & 3.125 & 0 & 0 & 2.5 & 0 \\
    6 & 1 & 8 & 0 & 0 & 2 & 0 & 0 \\ 
    \hline
    \end{tabular}
    \footnotetext{Here $p\up{D,e}_t=\sum_{v \in \mathcal{V}}p\up{D,e}_{vt}$ and $e\up{e}_t=\sum_{v \in \mathcal{V}}e\up{e}_{vt}$. Loads and market prices are the same as in Table~\ref{tab:res_sw_low}.}
    \label{tab:res_sw_low_dis}
\end{table}
\begin{table}
    \caption{Update of the inter-storage system after each market clearing for the second example of Sect.~\ref{sec:sw}, with discount}
    \centering
    \begin{tabular}{l|rr|rr}
    \hline
    \multicolumn{1}{l|}{} & \multicolumn{2}{c|}{$v=1$} & \multicolumn{2}{c}{$v=2$}\\ 
    \hline
    \multicolumn{1}{l|}{MC} & \multicolumn{1}{l}{$E_{v}\up{init}$(MWh)} & \multicolumn{1}{l|}{$S_v$ (\euro{}/MWh)} & \multicolumn{1}{l}{$E_{v}\up{init}$(MWh)} & \multicolumn{1}{l}{$S_v$ (\euro{}/MWh)}\\ 
    \hline
    1 & 2.5 & 31.25 &  &  \\
    2 & 2.5 & 20.31 &  &  \\
    3 & 2.5 & 13.20 &  &  \\
    4 & - & - &  &  \\
    5 & - & - & 2.5 & 1.56 \\
    6 & - & - & - &  -\\ 
    \hline
    \end{tabular}
    \footnotetext{Dashes indicate depletion.}
    \label{tab:res_upd3}
\end{table}

\subsection{Limitation} \label{sec:limit}
We have seen that the introduced MCP with VLBs performs better than the standard MCP in terms of ensuring cost recovery for the market participants. However, not all the limitations of the standard MCP are overcome, which we see here in a last illustrative example. 

We again consider a storage system with capacity $\overline{E}= 2.5$~MWh and initially empty. We have $\overline{P}\up{C}$ = $\overline{P}\up{D} = 3.5$~MW. In this example, $\eta\up{C} = \eta\up{D} = 1.0$, in order to facilitate the interpretation of the results. The rest of the data and the results are shown in Table~\ref{tab:res_eff}. There is one load and one generator. The market is cleared for two market clearings of three hours each. The final storage level is set to its optimal value, obtained from the ideal MCP. Note that the final level at the end of the second market clearing is equal to 2.5~MWh because of subsequent market clearings, which are not considered here.
\begin{sidewaystable}
    \caption{Data and results for the example of Sect.~\ref{sec:limit}}
    \centering
    \begin{tabular}{lr|rrrr|rrrrrr|rrrr|rrrrrr}
    \hline
    \multicolumn{1}{l}{} & \multicolumn{1}{l|}{} & \multicolumn{4}{c|}{Data} & \multicolumn{6}{c|}{Ideal MCP} & \multicolumn{4}{c|}{Standard MCP} & \multicolumn{6}{c}{MCP with VLBs}\\ 
    \hline
    \multicolumn{1}{l}{MC} & \multicolumn{1}{l|}{$t$} & \multicolumn{1}{l}{$D_{t}$} & \multicolumn{1}{l}{$U_{t}$} & \multicolumn{1}{l}{$P_{t}$} & $C_{t}$ & \multicolumn{1}{c}{$d_{t}$} & \multicolumn{1}{c}{$p_{t}$} & \multicolumn{1}{c}{$p\up{C}_{t}$} & \multicolumn{1}{c}{$p\up{D}_{t}$} & $e_t$ & \multicolumn{1}{l|}{$\lambda_t$} & $p\up{C}_t$ & \multicolumn{1}{c}{$p\up{D}_{t}$}  & $e_t$ & $\lambda_t$ & \multicolumn{1}{l}{$p\up{C,a}_t$} & \multicolumn{1}{l}{$p\up{D,a}_t$} & \multicolumn{1}{l}{$p\up{D,e}_t$} & \multicolumn{1}{l}{$e\up{a}_t$} & \multicolumn{1}{l}{$e\up{e}_t$} & \multicolumn{1}{l}{$\lambda_t$} \\ 
    \hline
    1 & 1 & 0 & 0 & 1 & 0 & 0 & 1 & 1 & 0 & 1 & 2 & 1 & 0 & 1 & 2 & 1 & 0 & 0 & 1 & 0 & 2 \\
    1 & 2 & 1 & 5 & 2.5 & 2 & 1 & 2.5 & 1.5 & 0 & 2.5 & 2 & 1.5 & 0 & 2.5 & 2 & 1.5 & 0 & 0 & 2.5 & 0 & 2 \\
    1 & 3 & 2 & 6 & 0 & 0 & 2 & 0 & 0 & 2 & 0.5 & [3,4] & 0 & 2 & 0.5 & [3,\textbf{6}] & 0 & 2 & 0 & 0.5 & 0 & [3,\textbf{6}] \\ 
    \hline
    2 & 1 & 2.5 & 4 & 2 & 3 & 2.5 & 2 & 0 & 0.5 & 0 & [3,4] & 0 & 0.5 & 0 & [3,4] & 0 & 0.5 & 0 & -0.5 & 0.5 & [3,4] \\
    2 & 2 & 0 & 0 & 1 & 2 & 0 & 1 & 1 & 0 & 1 & [3,7] & 1 & 0 & 1   & [3,7] & 1 & 0 & 0 & 0.5 & 0.5 & [3,7] \\
    2 & 3 & 4 & 7 & 5.5 & 3 & 4 & 5.5 & 1.5 & 0 & 2.5 & [3,7] & 1.5  & 0 & 2.5 & [3,7] & 1.5 & 0 & 0 & 2 & 0.5 & [3,7] \\
    \hline
    \end{tabular}
    \footnotetext{Here $p\up{D,e}_t=\sum_{v \in \mathcal{V}}p\up{D,e}_{vt}$ and $e\up{e}_t=\sum_{v \in \mathcal{V}}e\up{e}_{vt}$. The state of energy is in~MWh, prices in~\euro{}/MWh, and the rest is in~MW. Load values $d_{t}$ and generation values $p_{t}$ are the same in all the cases and are not repeated.}
    \label{tab:res_eff}
\end{sidewaystable}

We assume perfect foresight, and we can see that all quantities agree. As a consequence, social welfare is the same for all these MCPs. However, we observe that the prices in the last hour of the first market clearing can potentially be higher than in the ideal MCP for both other approaches. This difference is due to future market clearings not being considered for pricing in these cases. With the ideal MCP, the price for that hour is at most 4~\euro{}/MWh, which corresponds to the utility of the load in the first hour of the next market clearing, which the other approaches do not take into consideration.

We argue for a limited impact of future uncertain information on the formation of current prices, to ensure transparency. However, we see here that it comes with a cost that somebody will ultimately have to pay. The study of this trade-off is a topic for future research.

\section{Discussion and conclusion} \label{sec:5}

\noindent We introduced a novel procedure for clearing an energy market with non-merchant storage, using virtual linking bids. This is based on an artificial representation of the storage system, dividing it into a component for local arbitrage, within the current market clearing, and a component for arbitrage between market clearings. 

We showed that it outperforms traditional approaches when it comes to cost recovery.
Indeed, it ensures cost recovery for the storage system, even over multiple market clearings, which the standard MCP does not. More importantly, we showed that this property also stands when forecast errors are made when calculating the final state of energy of the storage system, which is usually the case. 

It still remains to study how and by whom this final state of energy should be determined. This should also include a study of the impact of the storage level on pricing. In particular, a critical next step would be to investigate the potential impact of a strategic choice of this level on prices and social welfare, and how it compares to having merchant storage.

We also discussed the impacts of uncertainty on social welfare compared to traditional approaches. In the case of forcing the discharge of the storage system at a disadvantageous price, we showed that our approach comes closer to closing the gap with an ideal MCP.

We also illustrated that this method does not solve the problem of accounting for prices in future market clearings, potentially leading to higher prices compared to an ideal MCP. This problem needs to be further studied.
Another limitation of the method is that it might happen that a stored quantity is kept for too long, due to a very high value. We showed that a discount factor could be applied to the value of stored energy over time to avoid this.

We used illustrative small-scale examples to give better intuition on how the method introduced behaves compared to the standard MCP. Since this method does not introduce non-linearities, the computational complexity for real instances is expected to be similar to traditional approaches. It only involves solving one more linear program for the update of storage values, which should also scale well.
This procedure was introduced on an idealized representation of a storage system and of the market clearing, and the promising results are a good motivation for extending it to more general setups, in particular to consider the degradation of the storage system and the provision of different services.

\section*{Declarations}
\paragraph{Competing interests} The authors have no competing interests to declare that are relevant to the content of this article.

\bibliography{ref}

\end{document}